\documentclass[12pt,leqno]{article}
\title{Some congruences on prime factors of class number of  extensions $K/\Q$, \  version 2.1}
\author{Roland Qu\^eme}

\usepackage{amsmath,amsthm,amstext,amsbsy,amsfonts}
\newtheorem{thm}{Theorem}[section]
\newtheorem{cor}[thm]{Corollary}

\newcommand{\N}{\mathbb{N}}

\newcommand{\Q}{\mathbb{Q}}

\newcommand{\modu}{\ \mbox{mod}\ }
\newcommand{\be}{\begin{equation}}
\newcommand{\ee}{\end{equation}}
\date{2002 December 20}
\begin{document}
%%% ====================================================================
\abstract $ $
\begin{itemize}
\item
This paper is a contribution to the description of the odd prime factors of the class number of the number fields.
\item
An  example  of the results obtained  is:

{\it Let $K/\Q$ be an abelian  extension with $N=[K:\Q]>1,\quad N$ odd.
Let $h(K)$ be the class number of $K$. Suppose that $h(K)>1$.
Let $p$ be a prime dividing $h(K)$.
Let $r_p$ be the rank of the $p$-class group of $K$.
Then $p\times (p^{r_p}-1)$ and $N$ are not coprime.}
\item
We give also in this paper a connection with Geometry of Numbers point of view. With  an explicit  geometric  upper bound $H_F$ of the class number $h(F)$ for any  field $F$, which is  given in this paper:

{\it Let $K/\Q$ be a Galois  extension with $[K:\Q]=N$. Let $h(K)$ be the class number of $K$. Suppose that $h(K)>1$. Suppose that $N$ has odd prime divisors. Let $n$ be an odd prime divisor of $N$.
Then  there exists a cyclic extension $K/F$ with  $n=[K:F]$.
Suppose that  $p>H_F$ is  a prime dividing  $h(K)$.
Let $r_p$ be the rank of the $p$-class group of $K$.
Then  $p\times (p^{r_p}-1)\equiv 0\modu n$.}
\item
The proofs are {\bf  elementary}.
We give  several verifications of results obtained for cyclic and abelian  extensions from the tables in  Washington\cite{was}, Schoof \cite{sch}, Masley \cite{mas}, Girstmair \cite{gir}, Jeannin \cite{jea} and the tables of cubic totally real number fields of the $ftp$ server
 {\it megrez.math.u-bordeaux.fr}.
\endabstract
\end{itemize}
%%% ====================================================================
\maketitle
%
%RRRRRRRR 10
%%% ====================================================================
\maketitle
%%% ====================================================================
%
%RRRRRRRR 10
%%% ====================================================================
\section{ On prime factors of class number of abelian extensions $K/\Q$}
The so called {\it rank theorem}, see Masley \cite{mas}, corollary 2.15 p. 305 is:

Suppose $M/P$ is a cyclic extension of degree $m$. Let $p$ be a prime which does not divide $h(E)$, class number of $E$  for
any field $E$ with $P\subset E\subset M,\quad E\not=M$, and which does not divide $m$. If $p|h(M)$ class number of $M$ then the rank $r_p$ of $p$-class group of $M$ is a multiple of $f$, the order of $p\modu m$.

In this subsection we use this theorem to get some congruences on prime factors of class number of abelian extensions $K/\Q$.

\subsection{ Some definitions}
\begin{itemize}
\item
Let $K/\Q$ be an abelian  extension  with $[K:\Q]=N$,  where $N=2^{\alpha_0}\times N_1,\quad N_1>1$  odd.
\item
Let us note in the sequel
the prime  decomposition
\begin{displaymath}
N=2^{\alpha_0}\times n_1^{\alpha_1}\times n_2^{\alpha_2}\times\dots\times n_k^{\alpha_k}.
\end{displaymath}
\item
Let  $h(K)$ be the class number of the field $K$. In
this paper, we  are studying, for fields $K$ with $h(K)>1$,  some congruences on the primes dividing $h(K)$.
\item
The extension $K/\Q$ being abelian, for any   prime $n$ dividing $N_1$,  there exists  at least one subfield $F$ of $K$ with $[K:F]=n$,  where $F/\Q$ is abelian and  $K/F$ is  cyclic.
Let us note $h(F)$  the class number of $F$.
\item
Let $p$ be a   prime dividing $h(K)$.
Let $r_p>0$ be the rank of the $p$-class group of $K$.
\end{itemize}
%%% ====================================================================
%
%RRRRRRRR 10
%%% ====================================================================
\subsection{Some results}
In  this subsection, we  give some explicit congruences on prime factors $p$ of class number $h(K)$ of abelian extensions $K/\Q$ in the general case and in the particular case where $N=[K:\Q]$ is odd.
\begin{thm}\label{t211291}
Let $K/\Q$ be an abelian extension with $[K:\Q]=2^{\alpha_0}\times N_1$, where $N_1>1$ is odd.
Let $h(K)$ be the class number of $K$. Suppose that $h(K)>1$.
Let $p$ be a prime dividing $h(K)$.
Let $r_p$ be the rank of the $p$-class group of $K$.
If $p\times (p^{r_p}-1)$ and $N_1$ are coprime,  then $p$ divides the class number $h(L)$, for  $L$ subfield of  $K$ with $[L:\Q]=2^{\alpha_0}$.
\begin{proof}$ $
\begin{itemize}
\item
There exists at least one cyclic extension $K/F$ of prime degree $n|N_1$.
\item
Suppose at first that $p\not|h(F)$:
Then hypotheses of   rank theorem are verified, see Masley \cite{mas} corollary 2.15 p. 305 :
\begin{itemize}
\item
From hypothesis,  $p$ does not divide $N_1$, so $p$ does not divide $n=[K:F]$.
\item
The extension $K/F$ is cyclic with
$ h(K)\equiv 0\modu p,\quad  h(F)\not\equiv 0\modu p$.
\item
There is no field $E$, different of $K$ and of $F$ with $F\subset E\subset K$.
\end{itemize}
From rank theorem, if $f$ is the order of $p \modu n$ then $f|r_p$ and so $p^{r_p}\equiv 1\modu n$,
and we are done.
\item
Suppose now that  $p$ divides $h(F)$ : then
$N^\prime=\frac{N}{n}=[F:\Q]$. We can pursue the same algorithm with abelian extension $F/\Q$ in place of
abelian  extension $K/\Q$
and $N^\prime$ in place of $N$,  up to find a prime divisor
$n_i,\quad 1\leq i\leq k,$  of $N_1$ dividing $p^{r_p}-1$,  or to get a subfield   $L$ of $K$ with
$[L:\Q]=2^{\alpha_0}$ and $p |h(L)$ , which achieves the proof.
\end{itemize}
\end{proof}
\end{thm}
%%% ====================================================================
%
%RRRRRRRR 10
%%% ====================================================================
When $N$ is odd, this leads to the particularly straightforward formulation:
\begin{cor}\label{c211301}
Let $K/\Q$ be an abelian  extension with $N=[K:\Q]>1,\quad N$ odd.
Let $h(K)$ be the class number of $K$. Suppose that $h(K)>1$.
Let  $p$ be a prime dividing $h(K)$.
Let $r_p$ be the rank of the $p$-class group of $K$.
Then $p\times (p^{r_p}-1)$ and $N$ are not coprime.
\end{cor}
%%% ====================================================================
%
%RRRRRRRR 10
%%% ====================================================================
\section{Geometry of Numbers point of view}
\begin{itemize}
\item
The next result connects rank theorem to Geometry of Numbers point of view.
Let $K/\Q$ be an algebraic extension
(here,  we now  {\bf don't} suppose that $K/\Q$ is abelian).
Suppose only that there exists a subfield $F$ of $K$ such that $K/F$ is a cyclic extension, with $[K:\Q]=N,\quad [K:F]=n$ where $n$ is an odd prime.
\item
Let   $D_F$ be the discriminant of $F$. Let $m=\frac{N}{n}$.
Then  $h(F)\leq \frac{2^{m-1}}{(m-1)!}\times \sqrt{|D_F|}\times (log(|D_F|))^{m-1}$,
see Bordell\`es \cite{bor} theorem 5.3 p. 4.
Let us note
\begin{equation}\label{e212141}
H_F=\frac{2^{m-1}}{(m-1)!}\times \sqrt{|D_F|}\times (log(|D_F|))^{m-1}.
\end{equation}
If a prime $p$ verifies
$p>H_F$ then $p\not| h(F)$.
\end{itemize}
%%% ====================================================================
%
%RRRRRRRR 10
%%% ====================================================================
\begin{thm}\label{t212141}
Let $K/\Q$ be an algebraic extension. Let $h(K)$ be the class number of $K$. Suppose that $h(K)>1$.
Suppose that there exists a cyclic extension $K/F$, where $n=[K:F]$ is an odd prime.
Let $H_F$ be a geometric upper bound of class number of $F$ given by relation (\ref{e212141}).
Suppose that  $p>H_F$ is  a prime dividing  $h(K)$.
Let $r_p$ be the rank of the $p$-class group of $K$.
Then  $p\times (p^{r_p}-1)\equiv 0\modu n$.
\begin{proof}
We apply rank theorem, see Masley\cite{mas} corollary 2.15 p 305, and upper bound
$H_F$ of class number $h(F)$ of field $F$ given in relation (\ref{e212141}).
\end{proof}
\end{thm}
%%% ====================================================================
%
%RRRRRRRR 10
%%% ====================================================================
\begin{cor}\label{c212151}
Let $K/\Q$ be a Galois  extension with $[K:\Q]=N$. Let $h(K)$ be the class number of $K$. Suppose that $h(K)>1$. Suppose that $N$ has odd prime divisors. Let $n$ be an odd prime divisor of $N$.
Then  there exists a cyclic extension $K/F$ with  $n=[K:F]$.
Let $H_F$ be a geometric upper bound of class number of $F$ given by relation (\ref{e212141}).
Suppose that  $p>H_F$ is  a prime dividing  $h(K)$.
Let $r_p$ be the rank of the $p$-class group of $K$.
Then  $p\times (p^{r_p}-1)\equiv 0\modu n$.
\begin{proof}
Immediate consequence of theorem \ref{t212141} and of Galois theory.
\end{proof}
\end{cor}
%%% ====================================================================
%
%RRRRRRRR 10
%%% ====================================================================
%%% ====================================================================
%
%RRRRRRRR 10
%%% ====================================================================
\section{Numerical examples}\label{ss22043}
The examples found to  check theses results are taken from:
\begin{itemize}
\item
the table of relative class numbers of cyclotomic number fields in Washington, \cite{was} p 412,  with some elementary MAPLE computations,
\item
the table of relative class number of cyclotomic number fields in Schoof, \cite{sch}
\item
the table of maximal real subfields $\Q(\zeta_l+\zeta_l^{-1})$ of $\Q(\zeta_l)$ for $l$ prime in Washington, \cite{was} p 420.
\item
the tables of relative class number of imaginary cyclic fields of Girstmair of degree 4,6,8,10. \cite{gir}.
\item
the tables of quintic number fields computed by Jeannin, \cite{jea}
\item
the tables of cubic totally real cyclic number fields of the Bordeaux University in the Server

{\it megrez.math.u-bordeaux.fr}.
\end{itemize}

All the results examined in these tables are in accordance with our theorems.
%%% ====================================================================
%
%RRRRRRRR 10
%%% ====================================================================
\subsection{Cyclotomic number fields $\Q(\zeta_u)$}

Let $u\in\N,\quad u>2$. Let $\zeta_u$ be a primitive $u^{th}$ root of unity.
Here we have $N=\phi(u)$, where $\phi$ is the Euler indicator.
The cyclotomic number fields $\Q(\zeta_u)$ of the examples are taken  with $2\| N=\phi(u)$,
except the example $\zeta_u$ with $u=572,\quad N=\phi(u)=2^4.3.5$.

\begin{itemize}
\item

$\Q(\zeta_u), \quad u=59 ,\quad N=\phi(u)=58=2.29 :$

$h^- = 3 . (2.29+1). (2^3 . 29 +1),\quad  h^+=1.$

$3 =h(\Q(\sqrt{-59}).$
\item

$\Q(\zeta_u), \quad u=71, \quad N=\phi(u)=70=2.5.7$ :

$h^- = 7^2.(2^3.5.7.283+1), \quad h^+=1. $

\item
$\Q(\zeta_u), \quad u=79,\quad  N=\phi(u)=78=2.3.13$ :

$ h^-= 5.(2^2.13+1)(2.3^2.5.13.17.19+1),\quad h^+=1.$

$5=h(\Q(\sqrt{-79})$.

\item
$\Q(\zeta_u), \quad u=83,\quad N=\phi(u)=82=2.41$ :

$ h^-=3.(2^2.41.1703693+1),\quad h^+=1$.

$3=h(\Q(\sqrt{-83}).$
\item
$\Q(\zeta_u),\quad u=103,\quad  N=\phi(u)=102=2.3.17$:

$h^-= 5.(2.3.17+1)(2^2.3.5.17+1)(2.3^2.5.17.11273+1)$,

$h^+=1$.

$5=h(\Q(\sqrt{-103}).$
\item
$\Q(\zeta_u), \quad u=107,\quad N=\phi(u)=106=2.53$ :

$h^-=(2.7.53+1)(2.3.31.53+1)(2^6.23.37.53+1), \quad h^+=1.$
\item
$\Q(\zeta_u), \quad u=121,\quad N=\phi(u)=110=2.5.11$ :

$h^-=(2.3.11+1)(2^5.11+1)(2^2.5.11.13+1)(2^2.3^2.5.11.13)$,

$h^+=1.$

\item
$\Q(\zeta_u), \quad u=127,\quad N=\phi(u)=126=2.3^2.7$ :

$ h^-=5.(2^2.3+1)(2.3.7+1)(2.3.7.13+1)(2.3^2.7^2+1)\times$

$(2.3^4.19+1)(2.3^2.7.4973+1), \quad h^+=1.$

$5=h(\Q(\sqrt{127})$.

\item
$\Q(\zeta_u), \quad u=131,\quad N=\phi(u)=130=2.5.13$ :

$ h^-=3^3.5^2.(2^2.13+1)(2.5.13+1)(2^2.5^2.13+1)\times$

$(2^2.5^2.13.29.151.821+1), \quad h^+=1$.

Observe that for $p=3$, we have  $p^3\|h$ and the group $C_p$ is not cyclic as it is seen in Schoof \cite {sch} table 4.2 p 1239 , where $r_p=3$. The rank theorem in that case shows that $f=r_p=3$ and so that $3^f=27\equiv 1 \modu 13$.

\item
$\Q(\zeta_u), \quad u=139,\quad N=\phi(u)=138=2.3.23$ :

$ h^-=3^2(2.23+1)(2^2.3.23+1)(2.3.7.23+1)(2^2.3.23.4307833+1),$

$h^+=1.$
\item
$\Q(\zeta_u), \quad u=151,\quad N=\phi(u)=150=2.3.5^2$ :

$ h^-=(2.3+1)(2.5+1)^2.(2^2.5.7+1)(2.3.5^2.173+1)\times$

$(2^2.3.5^4.7.23+1)(2^2.3.5^2.7.13.73.1571+1),\quad h^+=1.$

For $p=2.5+1=11$, from Schoof, see table 4.2 p 1239,  the group is not cyclic, $r_p=2$. From rank theorem, we can only assert that $f|r_p=2$, so that $f=1$ or $f=2$, so only that $p^2\equiv 1 \modu g=5$. We see in this numerical application that in that case  $f=1$ and $f<r_p$ strictly.

\item
$\Q(\zeta_u), \quad u=163,\quad N=\phi(u)=150=2.3^4$ :

$ h^-=(2^2.3^2.5+1)(2.3^4.11.13+1)(2^5.3^5.47+1)\times$

$(2.3^4.17.19.29.71.73.56179+1),\quad h^+=2\times 2$.
\item
$\Q(\zeta_u), \quad u=167,\quad N=\phi(u)=166=2.83$ :

$ h^-=11.(2.3.83+1)(2.83.22107011.1396054413416693+1),$

$h^+= 1.$

$11=h(\Q(\sqrt{167})$.
\item
$\Q(\zeta_u), \quad u=179,\quad N=\phi(u)=178=2.89$ :

$ h^-=(2^2.3.89+1)\times$

$(2^4.5.89.173.19207.155731.3924348446411+1),\quad h^+=1.$

\item $\Q(\zeta_u), \quad u=191,\quad N=\phi(u)=190=2.5.19$:

$ h^-= (2.5+1).13.(2.19^2.71+1)(2.3.5.19.277.3881+1),$

$h^+= (2.5+1)$.

$13=h(\Q(\sqrt{-191})$.

Here, we note that the prime $p=11$ corresponds to $11$-class group of
$\Q(\zeta_{191}+\zeta_{191}^{-1})$.

\item

$\Q(\zeta_u), \quad u=199,\quad N=\phi(u)=198=2.3^2.11$ :

$ h^-=3^4(2.3^2+1)(2.3.11^2+1)(2^2.3.11.23.8447)(2^4.3^2.11.13.17.331.1789),$

$h^+=1$.

\item
$\Q(\zeta_u),\quad u=572=2^2.11.13,\quad N=\phi(u)=2^4.3.5,$

$h^-=3.5^2.(2.3+1)(2.3^2+1)^2(2.3.5+1)(2^3.5+1)(2^2.3.5+1)^2$

$(2^2.3.5.7+1)(2^2.3.5.11+1)(2.3^2.5.307+1)(2.3^2.5.11.73+1)$

$(2^2.3^2.5^2.31+1)(2^2.3^5.5.7.53.263+1)$

We don't know $h^+$.

\end{itemize}
{\bf Remarks:}
\begin{itemize}
\item
Let $[K:\Q]=2^{\alpha_0}\times n_1^{\alpha_1}\times\dots\times n_k^{\alpha_k}$. We observe in  Washington  \cite{was}, tables  of relative class numbers  p. 412 and of real class numbers p. 421  that, when $p$ is large, $p-1$ is divisible by several  or all primes in the set
$\{n_i\ | \  i=1,\dots,k\}$. This observation complies with Geometry of Number corollary
\ref{c212151} p. \pageref{c212151}.
\item
Observe that frequently , we get $f=r_p=1$ and so $p\times (p-1)\equiv 1\modu n$. In our examples we get
only one example $\Q(\zeta_{131})$ with $p=3,\quad f=r_p=3$.
\end{itemize}
%%% ====================================================================
%
%RRRRRRRR 10
%%% ====================================================================
\subsection{Real class number $\Q(\zeta_l+\zeta^{-1}_l)$}
The examples are obtained from the table of real class number in Washington, \cite{was} p 420.
Here, $l$ is a prime, the class number  $h_\delta$ is the conjectured value of the class number $h^+$ of
$\Q(\zeta_l+\zeta_l^{-1})/\Q$ with a minor incertitude on an extra factor. But the factor $h_\delta$ must verify our theorems.
We extract some examples of the table with $2\| l-1$. Here we have $N=\frac{l-1}{2}$ and thus
$N\not \equiv 0 \modu 2$.
\begin{itemize}
\item
$l=191,\quad l-1=2.5.19,\quad h_\delta=(2.5+1)$

\item
$l=1063,\quad l-1=2.3^2.59,\quad h_\delta=(2.3+1).$

\item
$l=1231,\quad l-1=2.3.5.41,\quad h_\delta=(2.3.5.7+1).$

\item
$l=1459,\quad l-1=2.3^6,\quad h_\delta=(2.3.41+1).$

\item
$l=1567,\quad l-1=2.3^3.29,\quad h_\delta=(2.3+1).$

\item
$l=2659,\quad l-1=2.3.443,\quad h_\delta=(2.3^2+1).$

\item
$l=3547,\quad l-1=2.3^2.197,\quad h_\delta=((2.3^2+1)(2.3^2.7^2+1).$

\item
$l=8017,\quad l-1=2^4.3.167$,

$h_\delta=(2.3^2+1)(2.3^2.7^2+1)(2^2.3^3+1).$

\item
$l=8563,\quad l-1=2.3.1427,\quad h_\delta=(2.3+1)^2$. We have $r_p=1$ or $r_p=2$ so $f=1$ or $f=2$. We can conclude from theorem \ref{t211291} p. \pageref{t211291}, that $p^2\equiv 1\modu 3$.

\item
$l=9907,\quad l-1=2.3.13.127,\quad h_\delta=(2.3.5+1).$
\end{itemize}

%%% ====================================================================
%
%RRRRRRRR 10
%%% ====================================================================

\subsection {Cubic fields $K/\Q$ cyclic and totally real.}
Here, we have $N=3$. Note that in that case discriminants are square in $\N$.
\begin{itemize}
\item Discriminant $D_1=3969=(3^2.7)^2$, \quad$h=3$
\item Discriminant $D_2=3969=(3^2.7)^2$, \quad $h=7=2.3+1$
\item Discriminant $D_1=8281=(7.13)^2$, \quad $h=3$
\item Discriminant $D_1=13689=(3^2.13)^2$,\quad $h=3$
\item Discriminant $D_2=13689=(3^2.13)^2$, \quad $h=13=2^2.3+1$
\item Discriminant $D_1=17689=(7.19)^2$, \quad $h=3$
\end{itemize}
%\end{itemize}
%%% ====================================================================
%
%RRRRRRRR 10
%%% ====================================================================
\subsection{Totally real cyclic fields of prime conductor $<100$}
We have found few numeric results in the literature. We refer to Masley, \cite{mas}, Table 3 p 316.
In these results, $K/\Q$ is a real cyclic field with $[K:\Q]=N$, with conductor $f$, with root of discrimant $Rd$ and class number $h$.
\begin{itemize}
\item
$f=63,\quad N=3,\quad Rd=15.84,\quad \quad h=3$
\item
$f=63,\quad N=3,\quad Rd=15.84,\quad \quad h=3$
\item
$f=63,\quad N=6,\quad Rd=26.30,\quad \quad h=3$
\item
$f=63,\quad N=6,\quad Rd=26.30,\quad \quad h=3$
\item
$f=91,\quad N=3,\quad Rd=20.24,\quad \quad h=3$
\item
$f=91,\quad N=3,\quad Rd=20.24,\quad \quad h=3$
\item
$f=91,\quad N=6,\quad Rd=31.03,\quad \quad h=3$
\item
$f=91,\quad N=6,\quad Rd=31.03,\quad \quad h=3$
\end{itemize}
%%% ====================================================================
%
%RRRRRRRR 10
%%% ====================================================================
\subsection{Lehmer quintic cyclic field}
The prime divisors of the $82$ cyclic number fields of the table in Jeannnin \cite{jea}, with conductor $f<3000000$,  are, at a glance, of the form $p=2$ or $p=5$ or $p\equiv 1 \modu 10$, which clearly verifies corollary \ref{c211301} p.\pageref{c211301}.
%%% ====================================================================
%
%RRRRRRRR 10
%%% ====================================================================
\subsection{Decimic imaginary cyclic number fields with conductor between  $9000$ and $9500$}
This example is obtained from the tables of Girstmair, \cite{gir}. $f$ is a prime conductor, $h$ is the factorization of the class number $K/\Q$.
\begin{itemize}
\item
$f=9011$, $h=3.(2.5+1).(2.3.5.52201+1)$.

$3$ divides the class number of $\Q(\sqrt{-9011})$.
\item
$f=9081$, $h=3.7.(2^3.5+1)$.

$3,7$ divide the class number of $\Q(\sqrt{-9081})$.
\item
$f=9151$, $h=67.(2^2.5.1187+1)$.

$67$ divides the class number of $\Q(\sqrt{-9151})$.
\item
$f=9311$, $h=97.(2.5.5689+1)$.

$97$ divides the class number of $\Q(\sqrt{-9311})$.
\item
$f=9371$, $h=7^2.(2.3.5^2+1).(2.3^3.5+1)$.

$7$ divides the class number of $\Q(\sqrt{-9371})$.
\item
$f=9391$, $h=5^2.(2^4.3.5.7.71+1)$.

\item
$f=9431$, $h=7.13.(2.3.5+1)(2^4.3.5.71+1)$.

$7,13$ divide the class number of $\Q(\sqrt{-9431})$.
\item
$f=9491$, $h=3^2.5^2.(2.3.5+1)(2^2.3^2.5.13+1)$.

$3$ divides the class number of $\Q(\sqrt{-9491})$.
\end{itemize}
%%% ====================================================================
\maketitle
%
%RRRRRRRR 10
%%% ====================================================================
%%% ====================================================================
%
%RRRRRRRR 10
%%% ====================================================================
%
%RRRRRRRR 560
%%% ====================================================================

2002 december 20

***************

Roland Qu\^eme

13 avenue du ch\^ateau d'eau

31490 Brax

France

e-mail : roland.queme@free.fr

home page: http://roland.queme.free.fr/index.html


\begin{thebibliography}{9}
\bibitem{bor}, O. Bordell\`es, \textit{Explicit upper bounds for the average order of $d_n(m)$ and application to class number}, Journal of Inequalities in Pure and Applied Mathematics, vol 3, issue 3, article 38, 2002, http://jipam.vu.edu.au/
\bibitem{gir} K. Girstmair, \textit{The relative class number of imaginary cyclic fields of degree 4, 6, 8 and 10}, Math. Comp., 61, 204, 1993, p 881-887.
\bibitem{jea} S. Jeannin, \textit{Nombre de classes et unit\'es des corps de nombres cycliques quintiques d'E. Lehmer}, J. Th\'eorie des Nombres de Bordeaux, 1996, 8, pp75-92.
\bibitem{leh} D.H. Lehmer, \textit{Prime factors of cyclotomic class numbers}, Math. Comp., 31, 138, 1977, pp 599-607.
\bibitem{mas} J.M. Masley, \textit{Class number of real cyclic number fields with small conductor}, Compositio Mathematica, 37, 3., 1978, pp 297-319.
\bibitem{sch} R. Schoof, \textit{Minus class groups of the fields of $l-$th roots of unity}, Math. Comp., 67, 223, 1998, pp 1225-1245.
\bibitem{was} L.C. Washington, \textit{Introduction to cyclotomic field, second edition}, Springer, 1997.
\end{thebibliography}
\end{document}